\documentclass{amsart}
\usepackage{amscd}
\usepackage{a4wide}
\usepackage{amssymb,amsmath}  
\usepackage{latexsym} 

\DeclareFontEncoding{OT2}{}{} 
\newcommand{\textcyr}[1]{%
 {\fontencoding{OT2}\fontfamily{wncyr}\fontseries{m}\fontshape{n}\selectfont #1}}

\newcommand{\sha}{{\mbox{\textcyr{Sh}}}}
\newcommand{\q}{\mathfrak{q}}
\newcommand{\p}{\mathfrak{p}}

\newcommand{\Q}{\mathbf{Q}}
\newcommand{\Z}{\mathbf{Z}}
\newcommand{\F}{\mathbf{F}}

\newcommand{\Ps}{\mathbf{P}}

\newcommand{\ra}{\rightarrow}

\newcommand{\str}{\mathcal{O}}

\renewcommand{\bar}[1]{\overline{#1}}
    \newtheorem{Lem}{Lemma}[section]
    \newtheorem{Prop}[Lem]{Proposition}
    \newtheorem{Thm}[Lem]{Theorem}
    \newtheorem{Cor}[Lem]{Corollary}

   \theoremstyle{definition}
    \newtheorem{Def}[Lem]{Definition}
    
    \theoremstyle{remark}
    
    \newtheorem{Rem}[Lem]{Remark}
    \newtheorem{Not}[Lem]{Notation}

    \DeclareMathOperator{\Gal}{Gal}
\DeclareMathOperator{\rank}{rank}

\DeclareMathOperator{\pri}{prime}

\begin{document}

\title[Tate-Shafarevich groups]{The $p$-part of Tate-Shafarevich
groups of elliptic curves can be arbitrarily large}

\author{Remke Kloosterman}
\address{Department of Mathematics and Computer Science\\ University of Groningen\\ PO Box 800\\ 9700 AV  Groningen\\ The Netherlands}
\email{r.n.kloosterman@math.rug.nl}

\subjclass{Primary 11G05; Secondary 11G18}
\keywords{Tate-Shafarevich group, elliptic curve, abelian variety}

\thanks{The author wishes to thank Jaap Top for many useful conversations on this topic. The author wishes to thank Jasper Scholten for suggesting \cite{Mi}. The author wishes to thank Stephen Donnolley for pointing out a mistake in an earlier version.}

\date{\today}

\begin{abstract} In this paper it is shown that for every prime $p>5$ the dimension of the $p$-torsion in the Tate-Shafarevich group of $E/K$ can be arbitrarily large, where $E$ is an elliptic curve defined over a number field $K$, with $[K:\Q]$ bounded by a constant depending only on $p$. From this we deduce that the dimension of the $p$-torsion in the Tate-Shafarevich group of $A/\Q$ can be arbitrarily large, where $A$ is an abelian variety, with $\dim A$ bounded by a constant depending only on $p$. 
\end{abstract}

\maketitle
\markboth{KLOOSTERMAN}%
         {TATE-SHAFAREVICH GROUPS}

\renewcommand{\baselinestretch}{1.1}
\renewcommand{\arraystretch}{1.3}


\section{Introduction}
For the notations used in this introduction we refer to Section~\ref{Selm}.

The aim of this paper is to give a proof of

\begin{Thm}~\label{MainThm} There is a function $g: \Z \ra \Z$ such that for every prime number $p$  and every $k\in \Z_{>0}$ 
there exist infinitely many  pairs $(E,K)$, with $K$ a number field of degree at most $g(p)$ and $E/K$ an elliptic curve, such that
\[ \dim_{\F_p} \sha(E/K)[p]>k. \]
\end{Thm}

A direct consequence is:

\begin{Cor}\label{MainCor} For every prime number $p$ and every $k\in \Z_{>0}$ there exist infinitely many non-isomorphic abelian varieties $A/\Q$, with $\dim A \leq g(p)$ and $A$ is simple over $\Q$,  such that 
\[ \dim_{\F_p} \sha(A/\Q)[p]>k. \]
\end{Cor}

In fact, a rough estimate using the present proof reveals that $g(p)=\str(p^4)$, for $p \ra \infty$.

For $p\in \{2,3,5\}$, it is known that the group $\sha(E/\Q)[p]$ can be arbitrarily large. (See \cite{Bo}, \cite{Cas1}, \cite{Fi} and \cite{Kr}.) So we may assume that $p>5$, in fact, our proof only uses $p>3$.

The proof is based on combining the strategy used in \cite{Fi} to prove that $\dim_{\F_5} \sha(E/\Q)[5]$ can be arbitrarily large and the strategy used in \cite{Kl} to prove that $\dim_{\F_p} S^p(E/K)$ can be arbitrarily large, where $E$ and $K$ vary, but $[K:\Q]$ is bounded by a function depending on $p$ of type $\str(p)$. 

In \cite{Kl} the strategy was to find a field $K$, such that $[K:\Q]$ is small and
a point $P \in X_0(p)(K)$ such that $P$ reduces to one cusp for many primes $\p$  and reduces to the other cusp for very few primes $\p$. Then to $P$ we can associate an elliptic curve $E/K$ such that an application of a Theorem of Cassels (\cite{Cas2}) shows that $S^p(E/K)$ gets large.

The strategy of \cite{Fi} can be described as follows. Suppose $K$ is a field with class number 1. Suppose $E/K$ has a $K$-rational point of order $p$, with $p>3$ a prime number. Let $\varphi: E \ra E'$ be the isogeny obtained by dividing out the point of order $p$. Then one can define a matrix $M$, such that the $\varphi$-Selmer group is isomorphic to the kernel of multiplication on the left by $M$, while the $\hat{\varphi}$-Selmer group  is isomorphic to the kernel of multiplication on the right by $M$. One can then show that the rank of $E(K)$ and of $E'(K)$  is bounded by the number of split multiplicative primes minus twice the rank of $M$ minus 1.

Moreover, one can prove that if the matrix $M$ is far from being square, then the dimension of the $p$-Selmer group of one of the two isogenous curves is large. If one has an elliptic curve with two rational torsion points of order $p$ and $q$ respectively, one can hope that for one isogeny the associated matrix has high rank, while for the other isogeny the matrix is far from being square. Fisher uses points on $X(5)$ to find elliptic curves $E/\Q$ with two isogenies, one such that the associated matrix has large rank, and the other such that the 5-Selmer group is large.

We generalize this idea to number fields, without the class number 1 condition.
 To do this we express the Selmer group attached to the isogeny  as the kernel of the multiplication on the left by some matrix $M$. In general, the matrix for the dual isogeny turns out to be different from the transpose of $M$.

\begin{Rem} An element of $S^p(E/K)$ corresponds to a Galois extension $L$ of $K(E[p])$ of degree $p$ or $p^2$, satisfying certain local conditions and such that the Galois group of $L/K(E[p])$ interacts in some prescribed way with the Galois group of $K(E[p])/K$. The examples of elliptic curves with large Selmer and large Tate-Shafarevich groups in \cite{Fi}, \cite{Kl} and this paper have one thing in common, namely that the representation of the absolute Galois group of $K$ on $E[p]$ is reducible. It follows that the conditions on the interaction of the Galois group of $K(E[p])/K$ with the Galois group of $L/K(E[p])$ almost disappear. The next thing to do would be to produce examples of large $p$-Selmer groups, with $p>3$, with an irreducible Galois representation on $E[p]$. One could start with the case of elliptic curves with complex multiplication.\end{Rem}

The organization of this paper is as follows: In Section~\ref{Selm} we prove several lower and upper bounds for the size of $\varphi$-Selmer groups, where $\varphi$ is an isogeny with kernel generated by a rational point of prime order at least 5. In Section~\ref{modular} we use the modular curve $X(p)$ and the estimates from Section~\ref{Selm} to prove Theorem~\ref{MainThm}.

\section{Selmer groups}\label{Selm}

Suppose $K$ is a number field, $E/K$ is an elliptic curve and $\varphi: E \ra E'$ is an isogeny defined over $K$. Let $H^1(K,E[\varphi])$ be the first cohomology group of the Galois module $E[\varphi]$. 

\begin{Def} The $\varphi$-Selmer group of $E/K$ is
\[ S^\varphi (E/K) := \ker  H^1(K,E[\varphi]) \ra \prod_{\p \pri} H^1(K_\p,E).  \]
and the Tate-Shafarevich group of $E/K$ is
\[ \sha(E/K) := \ker H^1(K,E) \ra \prod_{\p \pri} H^1(K_\p,E). \]
\end{Def} 

\begin{Not} \label{nota} For the rest of this section fix a prime number $p>3$, a number field $K$ such that $\zeta_p \in K$ and an elliptic curve $E/K$ such that there is a non-trivial point $P\in E(K)$ of order $p$. Let $\varphi: E \ra E'$ be the isogeny obtained by dividing out $\langle P \rangle$. 

To $\varphi$ we associate three sets of primes. Let $S_1(\varphi)$ be the set of split multiplicative primes, not lying above $p$, such that $P$ is not in the kernel of reduction. Let $S_2(\varphi)$ be the set of split multiplicative primes, not lying above $p$, such that $P$ is in the kernel of reduction. Let $S_3(\varphi)$ be the set of primes above $p$.

Suppose $\mathcal{S}$ is a finite sets of finite primes. Let \[K(\mathcal{S},p):= \{ x \in K^*/K^{*p} | v_\p(x) \equiv 0 \bmod p \; \forall \p \not \in \mathcal{S} \}. \]
Let $H^1(K,M;\mathcal{S})$ the subgroup of $H^1(K,M)$ of all cocycles not ramified outside $\mathcal{S}$.

For any cocycle $\xi\in H^1(K,M)$ denote $\xi_\p :=res_\p(\xi) \in H^1(K_\p,M)$. Let $\delta_\p$ be the boundary map 
\[ E'(K_\p)/\varphi(E(K_\p)) \ra H^1(K_\p,E[\varphi]).\]

Let $C_K$ denote the class group of $K$.
\end{Not}

Note that $S_1(\hat{\varphi})=S_2(\varphi)$ and $S_2(\hat{\varphi})=S_1(\varphi)$.

\begin{Prop} \label{PropphiSelmdesccoh} We have
\[ S^\varphi(E/K) =\{ \xi \in H^1(K,E[\varphi]; S_1(\varphi) \cup S_3(\varphi)) \mid \xi_
\p=0 \; \forall \p \in S_2(\varphi)  \mbox{ and } \xi_\p \in \delta_\p(E'(K_\p)/\varphi(E(K_\p)))\; \forall \p \in S_3(\varphi)\}.\]
\end{Prop}

\begin{proof} Note that if $p$ divides the Tamagawa number $c_{E,\p}$ then the reduction at $\p$ is split multiplicative. If this is the case then $c_{E,\p}/c_{E',\p}\neq 1$. This combined with $\dim H^1(K_\p,E[\varphi]) \leq 2$ (for $\p \nmid (p)$) and \cite[Lemma 3.8]{Sch} gives that $\iota_\p^*:H^1(K_\p,E[\varphi])\ra H^1(K_\p,E)$ is either injective or the zero-map. A closer inspection of \cite[Lemma 3.8]{Sch} combined with \cite[Proposition 3]{Kl} shows that $\iota_\p^*$ is injective if and only if $\p \in S_2(\varphi)$.

The Proposition follows then from \cite[Proposition 4.6]{SS}. \end{proof}

\begin{Prop} \label{PropphiSelmdesc} We have
\[ S^\varphi(E/K) \subset \{ x \in K(S_1(\varphi)\cup S_3(\varphi),p) | x \in K_\p^{*p} \;  \forall \p \in S_2(\varphi) \} \]
and
\[  S^\varphi(E/K) \supset \{ x \in K(S_1(\varphi), p) | x \in K_\p^{*p} \;  \forall \p \in S_2(\varphi)\cup S_3(\varphi) \}.\]\end{Prop}

\begin{proof} This follows from the identification $E[\varphi]\cong\Z/p\Z \cong \mu_\p$, the fact $H^1(L,\mu_p) \cong L^*/L^{*p}$ for any field $L$ of characteristic different from $p$, and Proposition~\ref{PropphiSelmdesccoh}.
\end{proof}

\begin{Def} Let $\mathcal{S}_1$ and $\mathcal{S}_2$ be two disjoint finite sets of finite primes of $K$, such that none of the primes in these sets divides $(p)$. 

Let 
\[ T: K(\mathcal{S}_1,p) \ra \oplus_{\p \in \mathcal{S}_2} \mathcal{O}_\p^* / \mathcal{O}_\p^{*p} \]
be the $\F_p$-linear map induced by inclusion.


Let $m(\mathcal{S}_1,\mathcal{S}_2)$ be the rank of $T$.

In the special case of an isogeny $\varphi: E\ra E'$ with associated sets $S_1(\varphi)$ and $S_2(\varphi)$ as above we write 
$m(\varphi):=m(S_1(\varphi),S_2(\varphi))$.
\end{Def}

Note that $K$ does not admit any real embedding. So 
\[ \dim K(S,p) = \frac{1}{2} [K:\Q] + \# S+\dim_{\F_p} C_K[p]\]
for any set of finite primes $S$.

\begin{Prop}\label{PropSphibnds} We have
\[ \# S_1(\varphi)-\# S_2(\varphi) + \dim_{\F_p}  C_K[p]-\frac{1}{2} [K:\Q]  \leq \dim S^\varphi(E/K) \leq \# S_1(\varphi)+ \dim_{\F_p}  C_K[p]-m(\varphi) + \frac{3}{2} [K:\Q]. \]
\end{Prop}

\begin{proof} Using Proposition~\ref{PropphiSelmdesc} twice we obtain
\begin{eqnarray*} 
- \frac{1}{2} [K:\Q] + \# S_1(\varphi)+\dim_{\F_p} C_K[p] - \# S_2(\varphi) &\leq & \dim K(S_1(\varphi),p)-\# S_2(\varphi) - \# S_3(\varphi) \\
& \leq & \dim S^\varphi(E/K)\\
& \leq & \dim K(S_1(\varphi) \cup S_3(\varphi),p) - m(\varphi)\\
& \leq & \# S_1(\varphi) + \# S_3(\varphi) +\dim_{\F_p} C_K[p]-m(\varphi)+\frac{1}{2} [K:\Q]\\
& \leq & \# S_1(\varphi) +\dim_{\F_p} C_K[p]-m(\varphi)+\frac{3}{2} [K:\Q].\end{eqnarray*}	
\end{proof}

\begin{Lem} \label{LemPSelmUp} 
We have
\[ \rank E(K) \leq\# S_1(\varphi)+ \# S_2(\varphi)+2\dim C_K[p] +3 [K:\Q] - m(\varphi) - m(\hat{\varphi})-1. \]
\end{Lem}

\begin{proof} Note
\begin{eqnarray*} 1+\rank E(K) &\leq&\dim_{\F_p}  E(K)/pE(K) \leq \dim S^p(E/K)\leq \dim S^\varphi(E/K) + \dim S^{\hat{\varphi}}(E'/K)  \\ 
&\leq& \# S_1(\varphi)+ \#S_1(\hat{\varphi}) +2 \dim C_K[p] + 3[K:\Q] - m(\varphi) - m(\hat{\varphi}),
\end{eqnarray*} which proves the Lemma.
\end{proof}

By a theorem of Cassels we can compute the difference of $ \dim S^\varphi(E/K)$ and $\dim S^{\hat{\varphi}}(E'/K)$. We do not need  the precise difference, but only an estimate, namely
\begin{Lem} \label{LemCassels} There is an integer $t$, with $|t| \leq 2[K:\Q]+1$ such that
\[ \dim S^{\hat{\varphi}}(E'/K) = \dim S^\varphi(E/K)  - \#S_1(\varphi)+\# S_2(\varphi) + t. \]
\end{Lem}

\begin{proof} This follows from \cite{Cas2} and \cite[Proposition 3]{Kl}. 
\end{proof}

\begin{Lem} \label{LemSumLow}
\[ \dim S^{\varphi}(E/K) + \dim S^{\hat{\varphi}}(E'/K)  \geq | \# S_1(\varphi) - \# S_2(\varphi) | + 2 \dim_{\F_p} C_K[p] - 3[K:\Q]-1.\]
\end{Lem}

\begin{proof} We may assume that $\#S_1 \geq  \# S_2$. From Proposition~\ref{PropSphibnds} we know 
\[ \dim S^\varphi(E/K) \geq  \# S_1(\varphi) - \# S_2(\varphi) + \dim_{\F_p} C_K[p] -\frac{1}{2} [K:\Q]. \]
From this inequality and Lemma~\ref{LemCassels}   we obtain that
\[ \dim S^{\hat{\varphi}}(E'/K)\geq  \dim S^\varphi(E/K) - 2[K:\Q]-1 - \# S_1(\varphi) +\# S_2(\varphi)  \geq \dim_{\F_p} C_K[p] -\frac{5}{2}[K:\Q]-1 .\]
Summing both inequalities gives the Lemma.
\end{proof}

\begin{Lem} \label{LemSumisSha} Suppose $\dim S^\varphi(E/K) + \dim S^{\hat{\varphi}}(E/'K)=k+1$ and $\rank  E(K)=r$, then 
\[ \max(\dim \sha(E/K)[p],\dim \sha(E'/K)[p]) \geq \frac{(k-r)}{2}.\]
\end{Lem}

\begin{proof} This follows from the exact sequence
\[ 0 \ra E'(K)[\hat{\varphi}]/\varphi(E(K)[p]) \ra S^\varphi(E/K) \ra S^p(E/K) \ra S^{\hat{\varphi}} (E'/K) \ra \sha(E'/K)[\hat{\varphi}]/\varphi(\sha(E/K)[p]). \]
(See \cite[Lemma 9.1]{SS}.)
\end{proof}

\begin{Lem} \label{LemMatrixisSha} Let $\psi: E_1 \ra E_2$ be some  isogeny obtained by dividing out a $K$-rational point of order $p$, with $E_1$  an elliptic curve in the $K$-isogeny network of $E$. Then
\[ \max(\dim \sha(E/K)[p],\dim \sha(E'/K)[p]) \geq  - \min(\# S_1(\varphi),\# S_2(\varphi)) - 3 [K:\Q]-1  +\frac{1}{2}( m(\psi) + m( \hat{\psi})).\]
\end{Lem}

\begin{proof} Use Lemma~\ref{LemPSelmUp} for the isogeny $\psi$ to obtain the bound for the rank of $E(K)$. Then combine this with Lemma~\ref{LemSumLow} and Lemma~\ref{LemSumisSha} and use that 
\[ \#S_1(\varphi)+\#S_2(\varphi)=\#S_1(\psi)+\#S_2(\psi). \]
\end{proof}

\section{Modular curves}\label{modular}
In this section we prove Theorem~\ref{MainThm}. The following result will be used in the proof of Theorem~\ref{MainThm}.

\begin{Thm}[{\cite[Theorem 10.4]{HR}}]\label{HRt} Let $f\in \Z[X]$ be a polynomial of
degree  at least 1.  Suppose that for every prime $\ell $, there exists a $y \in
\Z/\ell \Z$ such that $f(y) \not \equiv 0 \bmod \ell$. Then there exists a
constant $n$ depending on the degree of $f$ and the degree of its
irreducible factors such that there exist infinitely many primes $\ell$,
such that $f(\ell)$ has at most $n$ prime factors. Moreover, there exist $\delta>0, d\in \Z$, such that
\[ \# \{y\in \Z  \mid 0 \leq y \leq x \mbox{ and the number of prime factors of } f(y) \leq n \} \leq \delta \frac{x}{\log^{d} x } \left(1+\str\left(\frac{1}{\sqrt{\log(x)}}\right)\right) \]
as $x \ra \infty$.
 \end{Thm}

\begin{Not} Denote $X(p)/\Q$ the curve parameterizing triples $(E,Q_1,Q_2)$ with $\{Q_1,Q_2 \}$ an ordered basis for $E[p]$.

By the work of V\'elu (see \cite{Velu}, \cite[Chapter 4]{Fis}) there exists cusps of $X(p)$ defined over $\Q$. Fix a cusp $T_0 \in X_p(\Q)$. Let $R_1\in X_0(p)$ be the unramified cusp, let $R_2 \in X_0(p)$ be the ramified cusp.

Let $\tilde{\pi}_i : X(p)\ra X_0(p)$ be the morphisms obtained by mapping $(E,Q_1,Q_2)$ to $(E, \langle Q_i \rangle)$.

Let $\pi_i : X(p) \ra X_0(p)$ be either $\tilde{\pi}_i$ or $\tilde{\pi}_i$ composed with the Atkin-Lehner involution on $X_0(p$), such that $\pi_i(T_0)=R_1$.

Let $P\in X(p)$ a point, not a cusp. Denote $\varphi_{P,i}$ the morphism corresponding to $\pi_i(P)$.\end{Not}

 Note that the maps $\pi_i$ are defined over $\Q$.

\begin{Def} Let $T$ be a cusp of $X(p)$. We say that $T$ is of type $(\delta,\epsilon) \in\{1,2\}^2$ if $\pi_1(T)=R_\delta$ and $\pi_2(T)=R_\epsilon$. \end{Def}

Note that being of type $(\delta,\epsilon)$ is invariant under the action of the absolute Galois group of $\Q$.

For all points $P\in X(p)(\bar{\Q})$ and all primes $\p\nmid (p)$ of $\bar{\Q}$ such that $P \equiv T \bmod \p$ implies that  $\p \in S_\delta(\varphi_{P,1})$ and $\p \in S_\epsilon(\varphi_{P,2})$.

\begin{proof}[Proof of  Theorem~\ref{MainThm}] Fix a component $Y$ of $X(p)/\Q(\zeta_p)$. Fix a (possibly singular) model $C_1$ for  $Y$ in $\Ps^2$, such that the line $X=0$ intersects $C_1$ only in  cusps of type $(1,1)$ and no other point, all  $x$-coordinates of other the cusps are distinct and finite, 
if two $x$-coordinates of cusps are conjugate under $\Gal(\Q(\zeta_p)/\Q)$, then the type of these cusps coincide
and all $y$-coordinates of the cusps are finite. Let $C=\cup_{\sigma \in \Gal(\Q(\zeta_p)/\Q)}C_1^\sigma$. Then $C/\Q$ is a model for $X(p)$.

 Denote $h$ the defining polynomial of $C$.

Let $f_{\delta,\epsilon} \in \Z[X]$ be the radical polynomial with as roots all $x$-coordinates of the cusps of type $(\delta,\epsilon)$ of $X(p)$ and content 1. After a transformation of the form $x \mapsto cx$, we may assume that $f_{1,2}(0)f_{2,2}(0)=1$ and $f_{1,2}f_{2,2}\in \Z[X]$. Let $n$ denote the constant of Theorem~\ref{HRt} for the polynomial $f_{1,2}f_{2,2}$. Note that the discriminant of $f_{1,1}f_{1,2}f_{2,1}f_{2,2}$ is non-zero.

Let $\mathcal{B}$ consist of $p$, all primes $\ell$ dividing the leading coefficient or the discriminant of $f_{1,1} f_{1,2} f_{2,1} f_{2,2}$, all primes $\ell$  smaller then the degree of $f_{1,2}f_{2,2}$ and all primes dividing the leading coefficient of $res(h,f_{2,1},x)$.

Let $\mathcal{P}_2$ be the set of primes not in $\mathcal{B}$ such that every irreducible factor of $f_{2,1}(x) (x^p-1)\bmod \ell $ and every irreducible factor of $res(h,f_{2,1},x) \bmod \ell$  has degree 1. Note that by Frobenius' Theorem (\cite{Ste}) the set $\mathcal{P}_2$ is infinite. The condition mentioned here, implies that if we take a triple $(x_0,\ell,x_0)$ with $x_0 \in \Z$, the prime $\ell \in \mathcal{P}_2$ divides $f_{2,1}(x_0)$ and $y_0$ is a zero of $h(x_0,y)$ then every prime $\q$ of $\Q(\zeta_p,y_0)$ over $\ell$ satisfies $f(\q/\ell)=1$.

Fix $\mathcal{S}_1$ and $\mathcal{S}_2$ two sets of primes, such that 
\[m(\mathcal{S}_1,\mathcal{S}_2)> 2k+4(n+3)\deg(h) (p-1) +2,\]
 $S_1 \cap \mathcal{B}=\emptyset$  and $S_2 \subset  \mathcal{P}_2$. (The existence of such sets follows from Dirichlet's theorem on primes in arithmetic progression and the fact that $\ell \in \mathcal{S}_2$ implies $\ell \equiv 1 \bmod p$.)

\begin{Lem}\label{cond}There exists an $x_0\in \Z$ such that
\begin{itemize} 
\item $x_0 \equiv 0 \bmod \ell$, for all primes $\ell$ smaller then the degree of $f_{1,2}f_{2,2}$ and all $\ell$ dividing the leading coefficient of $f_{1,2}f_{2,2}$.
\item $x_0\equiv 0 \bmod \ell$, for all $\ell \in \mathcal{S}_1$,
\item $f_{2,1}(x_0) \equiv 0 \bmod \ell$, for all $\ell \in \mathcal{S}_2$,
\item $f_{1,2}(x_0)f_{2,2}(x_0)$ has at most $n$ prime divisors.
\item $h(x_0,y)$ is irreducible.
\end{itemize}
\end{Lem}

\begin{proof}
The existence of such an $x_0$ can be proven as follows. Take an $a\in \Z$ satisfying the above three congruence relations. Take $b$ to be the product of all primes mentioned in the above congruence relations. Define $\tilde{f}(Z)=f_{1,2}(a+bZ)f_{2,2}(a+bZ)$. We claim that the content of $\tilde{f}$ is one. Suppose $\ell$ divides this content. Then $\ell$ divides the leading coefficient of $\tilde{f}$. From this one deduces that $\ell$ divides $b$. We distinguish several cases:
\begin{itemize}
\item If $\ell \in \mathcal{S}_i$ then  $f_{i,2}(a)\equiv 0 \bmod \ell$ and $\ell$ does not divide the discriminant of the product of the $f_{\delta,\epsilon}$, so we have $\tilde{f}(0)\equiv f_{1,2}(a)f_{2,2}(a)\not \equiv 0 \bmod \ell$. 

\item If $\ell$ divides $b$ and is not in $\mathcal{S}_1 \cup \mathcal{S}_2$ then $\tilde{f}(0) \equiv f_{1,2}(0)f_{2,2}(0) \equiv 1 \bmod \ell$. 
\end{itemize}
So for all primes $\ell$ dividing $b$ we have that $\tilde{f}\not \equiv 0\bmod \ell$. This proves the claim on the content of $\tilde{f}$.

Suppose $\ell$ is a prime smaller then the degree of $\tilde{f}$,  then $\tilde{f}(0) \equiv 1 \bmod \ell$. 
If $\ell$ is different from these primes, then there is a coefficient of $\tilde{f}$ which is not divisible by $\ell$ and the degree of $\tilde{f}$ is smaller then $\ell$. So for every prime $\ell$ there is an $z_\ell\in \Z$ with $\tilde{f}(z_\ell) \not \equiv 0 \bmod \ell$. From this we deduce that we can apply Theorem~\ref{HRt}. The constant for $\tilde{f}$ depends only on the degree of the irreducible factors of $\tilde{f}$, hence equals $n$. The set 
\[ \{ x_1 \in \Z \mid \tilde{f}(x_1) \mbox{ has at most $n$ prime divisors} \} \]
is not a thin set. So
\[ \mathcal{H}:=\{ x_1 \in \Z \mid \tilde{f}(x_1) \mbox{ has at most $n$ prime divisors and $h(a+bx_1,y)$ is irreducible} \} \]
is not empty by Hilbert's Irreducibility Theorem \cite[Chapter 9]{Ser}.) Fix such an $x_1\in \mathcal{H}$. Let $x_0=a+bx_1$. This proves the claim on the existence of such an $x_0$.
\end{proof}


Fix an $x_0$ satisfying the conditions of Lemma~\ref{cond}. Adjoin a root $y_0$ of $h(x_0,y)$ to $\Q(\zeta_p)$. 
Denote the field $\Q(\zeta_p,y_0)$ by $K_1$. Let $P$ be the point on $X(p)(K_1)$ corresponding to $(x_0,y_0)$. Let $E/K_1$ be the elliptic curve corresponding to $P$. Let $K=K_1(\sqrt{c_4(E)})$. Then if $\q$ is a prime such that $E/K_\q$ has multiplicative reduction, then $E/K_\q$ has split multiplicative reduction.

For every prime $\p$ of $K$ over $\ell\in \mathcal{S}_1$ we have that $P \bmod \q$ is a cusp of type $(1,1)$.  Over every prime $\ell \in \mathcal{S}_2$ there exists a prime $\q$ such that $P \bmod \q$ is a cusp of type $(2,1)$. From our assumptions on $x_0$ it follows that $p$ does not divide $f(\q/\ell)$. Let $\mathcal{T}_1$ consists of the primes of $K$ lying over the primes in $\mathcal{S}_1$. Let $\mathcal{T}_2$ be the set of primes $\q$ such that $\q$ lies over a prime in $\mathcal{S}_2$ and $P \bmod \q$ is a cusp of type $(2,1)$.

Note that the set of primes of $K$ such that $P$ reduces to a cusp of type $(*,2)$ has at most $n[K:\Q]$ elements.

We have the following diagram
\[ \begin{array}{ccc}
\Q(\mathcal{S}_1,p)  & \ra & \oplus_{\ell \in \mathcal{S}_2} \Z_\ell^* / \Z_\ell^{*p} \\
\downarrow & & \downarrow \\
K(\mathcal{T}_1,p) & \ra & \oplus_{\q \in \mathcal{T}_2} \str_{K_\q}^* / \str_{K_\q}^{*p}. \end{array}\]
Since $p\nmid f(\q/\ell)$ for all $\ell \in \mathcal{S}_2$, the arrow in the right column is injective. This implies 
\[m(\varphi_{P,1} /K) \geq m(\mathcal{T}_1,\mathcal{T}_2) \geq m(\mathcal{S}_1,\mathcal{S}_2)=2k+4(n+3)\deg(h)(p-1)+2.\]

Since $S_2(\varphi_{p,2}/K) \leq [K:\Q] n$ and $[K:\Q]\leq 2(p-1) \deg(h)$ we obtain by Lemma~\ref{LemMatrixisSha} that for some $E'$ isogenous to $E$ we have
\begin{eqnarray*} \dim_{\F_p} \sha(E'/K)[p] &\geq& -\#S_2(\varphi_{P,2}) -3[K:\Q]-1+\frac{1}{2} m(\mathcal{S}_1,\mathcal{S}_2) \\
&\geq& (n+3)[K:\Q]-1+\frac{1}{2} m(\mathcal{S}_1,\mathcal{S}_2)=k.\end{eqnarray*}

Note that $\deg(h)$ can be bounded by a function of type $\str(p^3)$, hence $[K:\Q]$ can be bounded by a function of type $\str(p^4)$.
\end{proof}


To finish, we prove Corollary~\ref{MainCor}.

\begin{proof}[Proof of Corollary~\ref{MainCor}] Let $E/K$ be an elliptic curve such that $\dim \sha(E/K)[p]\geq k g(p)$ and $[K:\Q]\leq g(p)$. 

Let $R:=Res_{K/\Q}(E)$ be the Weil restriction of scalars of $E$. Then by \cite[Proof of Theorem 1]{Mi}
\[ \dim \sha(R/\Q)[p]=\dim \sha(E/K)[p].\]
 From this it follows that there is a factor $A$ of $R$, with $\dim \sha(A/\Q)[p] \geq k$.
\end{proof}

\end{document}